\documentclass[12pt]{article}
\usepackage{amsmath,enumerate} 
\usepackage{epsfig}
\usepackage{amssymb} 
\newtheorem{solution}{Exercise}
\newtheorem{lm}{{\bf Lemma}}[section]
\newtheorem{theor}[lm]
{{\bf Theorem}}
\newtheorem{deff}[lm]{{\bf Definition}}
\newtheorem{cor}[lm]{{\bf Corollary}}
\newtheorem{prop}[lm]{{\bf Proposition}}
\newtheorem{remark}[lm]{Remark}
\newtheorem{example}[lm]{{\bf Example}}

\newcommand{\lem}[1]{\begin{lm}$\!\!$\sep{\it #1}\end{lm}}
\newcommand{\theo}[1]{\begin{theor}$\!\!$\sep{\it #1}\end{theor}}
\newcommand{\coro}[1]{\begin{cor}$\!\!$\sep{\it #1}\end{cor}}

\newcommand{\propo}[1]{\begin{prop}$\!\!$\sep{\it #1}\end{prop}}
\newcommand{\rem}[1]{{\begin{remark}$\!\!$\sep{\rm #1}\end{remark}}}

\newcommand{\rf}[1]{{\rm(\ref{#1})}}

\newcommand{\pr}{\noindent {\sc Proof} --- \ }

\newcommand{\ep}{\hfill \framebox[2mm]{\ } \medskip}
\newcommand{\sep}{{\rm .} --- \ }

\newcommand{\be}{\begin{enumerate}}
\newcommand{\ee}{\end{enumerate}}

\newcommand{\bi}{\begin{itemize}}
\newcommand{\ei}{\end{itemize}}

\renewcommand{\leq}{\leqslant}
\renewcommand{\geq}{\geqslant}

\newcommand{\f}{\varphi}

\newcommand{\R}{\mathbb{R}}

\newcommand{\e}{\varepsilon}

\newcommand{\g}{{\bf g}}

\begin{document} 

\title{\bf Similarity solutions of mixed convection boundary-layer flows in a porous medium}
\author{M. A\"iboudi,  I. Bensari-Khelil, B. Brighi\,\footnote{\,Corresponding author}} 

\date{}
\maketitle

\begin{abstract}
The similarity differential equation $f'''+ff''+\beta f'(f'-1)=0$ with $\beta>0$ is considered.
This differential equation appears in the study of mixed convection boundary-layer flows over a vertical surface embedded in a porous medium. In order to prove the existence of solutions satisfying the boundary conditions $f(0)=a\geq0$, $f'(0)=b\geq0$ and $f'(+\infty)=0$ or $1$, we use shooting and consider the initial value problem consisting of the differential equation and the initial conditions $f(0)=a$, $f'(0)=b$ and 
$f''(0)=c$. 
For $0<\beta\leq1$, we prove that there exists a unique solution such that $f'(+\infty)=0$, and infinitely many solutions such that $f'(+\infty)=1$. For $\beta>1$, we give only partial results and show some differences with the previous case.
\end{abstract}

\section{Introduction}

Let $\beta\in\R$. 
We consider the third order autonomous nonlinear differential equation 
\begin{equation}
f'''+ff''+\beta f'(f'-1)=0. \label{eq}
\end{equation}
In fluid mechanics, in the study of mixed convection boundary-layer flows over a vertical surface embedded in a porous medium, such an equation arises in some situations where simplifying assumptions have been made ; see~\cite{aly}. Its solutions are called {\em similarity} solutions. 

\medskip

Equation~\rf{eq} is a particular case of the more general equation
\begin{equation}
f'''+ff''+\g(f')=0.\label{eqg}
\end{equation}
The most famous equation of this type is certainly the Blasius equation (see \cite{bla}), which corresponds to $\g=0$, and which has been extensively studied over the last hundred years ; see for example \cite{bfs} and the references therein.
   
For $\g(x)=\beta(x^2-1)$, this is the Falkner-Skan equation, introduced in 1931 for studying the boundary layer flow past a semi-infinite wedge, see the original paper \cite{falk} and \cite{hart} for a overview of mathematical results. 

For $\g(x)=\beta x^2$, this corresponds to free convection problems, see for example \cite{cheng} for the derivation of the model, and 
\cite{ab}, \cite{bbt}, \cite{bb}, \cite{asmq}, \cite{mmas}, \cite{bs}, \cite{bt}, \cite{g}, \cite{p}, \cite{tw} for different approaches of the mathematical analysis. 

The case where $\g(x)=\beta(x^2+1)$ is for the study of the boundary layer separation at a free stream-line, see \cite{a} and \cite{m}.
\medskip

Most of the time, these similarity equations are studied on the half line $[0,+\infty)$ and are associated to boundary conditions as $f(0)=a$, $f'(0)=b$ (or $f''(0)=c$) and a condition at infinity. This condition at infinity can be, either  $f'(t)\to\lambda$ as $t\to+\infty$, or $f'(t)\sim A\,t^{\nu}$ as $t\to+\infty$, where $A$ and $\nu$ are some positive constants,
or also $\vert f\vert$ is of polynomial growth at infinity. For more details, we refer to the introduction of \cite{rm} and to the references therein.

\medskip

The boundary value problems associated to the general equation~\rf{eqg}, with the condition that $f'$ tends to
$\lambda$ at infinity have been studied in~\cite{amuc} and in~\cite{rm}.
Let us notice that, if $\g(\lambda)\neq0$, then these boundary value problems do not have any solutions, and thus we must assume that $\g(\lambda)=0$ to have solutions. For example, in the case of mixed convection, i.e. $\g(x)=\beta x(x-1)$, the only relevant conditions are $f'(t)\to 0$ or $f'(t)\to1$ as $t\to+\infty$.
Results about existence, uniqueness and asymptotic behavior of {\em concave} or {\em convex} solutions to these boundary value problems are obtained, according to the sign of $\g$ between $b$ and $\lambda$.
Without further assumptions on $\g$, it is hopeless to have more precise results.  Nevertheless, the results of~\cite{rm} generalize the ones of~\cite{aml} and some of~\cite{guedda} about mixed convection problems.

\medskip

Let $a,b\in\R$ and $\lambda\in\{0,1\}$. 
We associate to equation \rf{eq}  the boundary value problem 
\begin{equation*}
\left \{\begin{array}
[c]{lll}
f'''+ff''+\beta f'(f'-1)=0\quad\mbox{on}\quad[0,+\infty)\\ \noalign{\vskip1,5mm}
f(0) = a\\ \noalign{\vskip1,5mm}
f'(0) =b\\ \noalign{\vskip1,5mm}
f'(t)\to\lambda\quad\mbox{as}\quad t\to+\infty
\end{array} \right.\eqno ({\cal P}_{\beta;a,b,\lambda})
\end{equation*}
Usually, the method to investigate such a boundary value problem is the shooting method, which consists of finding the values of a parameter $c$ for which the solution of \rf{eq} satisfying the initial conditions $f(0)=a$, $f'(0)=b$ and $f''(0)=c$, exists up to infinity and is such that $f'(t)\to\lambda$ as $t\to+\infty$. This approach is used in~\cite{aml} and~\cite{guedda}. 
In~\cite{aml}, the problem $({\cal P}_{\beta;a,b,1})$ is considered for $\beta<0$ and its is shown that this problem has a unique convex solution if $0<b<1$, and has a unique concave solution if $b>1$.
In~\cite{guedda}, for $\beta\in(0,1)$, $a=0$ and $b\in(0,\frac32)$, it is proven that the boundary value problem $({\cal P}_{\beta;a,b,1})$ has infinitely many solutions.

In~\cite{kyy},~\cite{y} and~\cite{yzd}, some results about the problem $({\cal P}_{\beta;a,b,1})$ are proven by introducing a singular integral equation obtained from \rf{eq} by a Crocco-type transformation.

\medskip

In the following, we will study the boundary value problems $({\cal P}_{\beta;a,b,0})$ and $({\cal P}_{\beta;a,b,1})$ for $\beta>0$, $a\geq0$ and $b\geq0$. 
In the case where $0<\beta\leq1$, we are able to get complete results (and so we improve the results of ~\cite{guedda}), while we only have partial results for $\beta>1$. On several occasions, we will use the results of~\cite{rm}, that sometimes we re-demonstrate, in our particular case, for the convenience of the reader.
 
The paper is organized as follows. In Section~\ref{s2}, general results about the solution of equation \rf{eq} are given. Section~\ref{s3} is devoted to the case where $b\geq 1$ and to the proofs of results that do not depend on whether $\beta\in(0,1]$ or $\beta>1$.
Section~\ref{s4} discusses in detail the case $\beta\in(0,1]$ and $b\geq 1$.
Section~\ref{s5} considers the case $\beta\in(0,1]$ and $0\leq b<1$, presents the results and how to prove them.  In Section~\ref{s6}, some results in the case $\beta>1$ are proven.

\section{\label{s2}Preliminary results}

To any $f$ solution of \rf{eq} on some interval $I$, we associate the function $H_f:I\to \R$ defined by
\begin{equation}
H_f=f''+f(f'-1)
.\label{hf}
\end{equation}
Then, we have
$H_f'=(1-\beta) f'(f'-1)
.$
\medskip

The following lemmas, concerning the solutions of the equation \rf{eq}, will be useful in the next sections. The proofs of some of them can be found in \cite{rm}.

\lem{\label{lc} Let $f$ be a solution of~\rf{eq} on some maximal interval $I$. If there exists $t_{0}\in I$ such that $f'(t_{0})\in\{0,1\}$ and  $f''(t_{0})=0$, then $I=\R$ and $f''(t)=0$ for all $t\in\R.$
}
\pr This follows immediatly from the uniqueness of solutions of initial value problem. Cf. \cite{rm}, Proposition 3.1, item 3.
\ep

\lem{\label{la} Let $\beta>0$ and $f$ be a solution of equation \rf{eq} on some interval $I$, such that $f'$ is not constant.
\vspace{-1mm}
\be[{\rm 1.}] \itemsep=0mm
\item If there exists $s<r\in I$ such that $f''(s)\leq0$ and 
$f'(f'-1)>0$ on $(s,r)$
then $f''(t)<0$ for all $t\in (s,r]$.

\item If there exists $s<r\in I$ such that $f''(s)\geq0$ and 
$f'(f'-1)<0$ on $(s,r)$
then $f''(t)>0$ for all $t\in (s,r]$.

\item If there exists $s<r\in I$ such that $f''<0$ on $(s,r)$ and 
$f''(r)=0$, then $f'(r)(f'(r)-1)<0$.

\item If there exists $s<r\in I$ such that $f''>0$ on $(s,r)$ and 
$f''(r)=0$, then $f'(r)(f'(r)-1)>0$.

\ee}
\pr Let $F$ denote any primitive function of $f$. From \rf{eq} we deduce the relation 
\begin{equation*}
(f''\exp F)'=-\beta f'(f'-1)\exp F.
\end{equation*}
All the assertions 1-4 follow easily from this relation and from Lemma~\ref{lc}. Let us verify the first and the third of these assertions.
For the first one, since $\psi=f''\exp F$ is decreasing on $[s,r]$, we have $f''(t)<f''(s)\exp(F(s)-F(t))\leq 0$ for all $t\in (s,r]$. For the third one, since $\psi<0$ on $(s,r)$ and $\psi(r)=0$, one has $\psi'(r)\geq 0$. This and Lemma~\ref{lc} imply that $f'(r)(f'(r)-1)<0$.
\ep

\lem{\label{lb} Let $f$ be a solution of~\rf{eq} on some maximal interval $(T_-,T_+)$. If $T_{+}$ is finite, then $f'$ and $f''$ are unbounded in any neighborhood of $T_+$.}
\pr Cf. \cite{rm}, Proposition 3.1, item 6.
\ep

\lem{\label{ld} Let $\beta\neq0$. If $f$ is a solution of~\rf{eq} on some interval $(\tau,+\infty)$ such that $f'(t)\to\lambda$ as $t\to+\infty$, then $\lambda\in\{0,1\}$. Moreover, if $f$ is of constant sign at infinity, then $f''(t)\to 0$ as $t\to+\infty$.
}
\pr Cf. \cite{rm}, Proposition 3.1, item 5 and 4. Let us notice that if $\lambda=1$, then $f$ is necessarily positive at infinity.
\ep

\lem{\label{le}  Let $\beta\neq0$. If $f$ is a solution of~\rf{eq} on some interval $(\tau,+\infty)$ such that  $f'(t)\to 0$ as $t\to+\infty$, then $f(t)$ does not tend to plus or minus infinity as $t\to+\infty$.}
\pr Assume for contradiction that $f(t)\to +\infty$ as $t\to+\infty$. Let $H=H_f$ be defined by \rf{hf}. Since $f'(t)\to 0$ as $t\to+\infty$, we deduce from the second assertion of Lemma~\ref{ld} that
$H(t)\sim-f(t)$ as $t\to+\infty$. This leads to a contradiction if $\beta=1$. If $\beta\neq1$, then we have
$H'(t)\sim(\beta-1) f'(t)$ as $t\to+\infty$, and hence $H(t)\sim(\beta-1) f(t)$ as $t\to+\infty$.
This is a contradiction, since $\beta\neq0$. The proof is the same if we assume that $f(t)\to -\infty$ as $t\to+\infty$. 
\ep

\lem{\label{lee} Let $\beta>0$ and $f$ be a solution of equation \rf{eq} on some right maximal interval $I=[\tau,T_+)$.
If $f\geq 0$ and $f'\geq0$ on $I$, then $T_+=+\infty$ and $\,f'$ is bounded on $I$.
}
\pr Let $L=L_f$ be the function defined on $I$ by
\begin{equation}
L(t)=3f''(t)^2+\beta (2f'(t)-3)f'(t)^2.\label{fonctionlf}
\end{equation}
Easily, using \rf{eq}, we obtain that
$L'(t)=-6f(t)f''(t)^2$ for all $t\in I$, and since $f\geq 0$ on $I$, this implies that $L$ is nonincreasing. Hence
$$\forall t\in I,\qquad \beta (2f'(t)-3)f'(t)^2\leq L(t)\leq L(\tau).$$
It follows that $f'$ is bounded on $I$ and, thanks to Lemma~\ref{lb}, that $T_+=+\infty$.
\ep

\lem{\label{pi} Let $\beta>0$ and $f$ be a solution of equation \rf{eq} on some right maximal interval $I=[\tau,T_+)$.
If $f(\tau)\geq 0$, $f'(\tau)\geq 1$ and $f''(\tau)>0$, then there exists $t_0\in(\tau,T_+)$ such that $f''>0$ on $[\tau,t_0)$ and $f''(t_0)=0$.
}
\pr Assume for contradiction that $f''>0$ on $I$. Then, $f'(t)\geq1$ and $f(t)\geq 0$ for all $t\in I$. We then have
\begin{equation}
f'''= -ff''-\beta f'(f'-1)\leq0. \label{a}
\end{equation}
It follows that $0<f''(t)\leq c$ for all $t\in I$ and hence, by Lemma~\ref{lb}, we have $T_+=+\infty$.
Next, let $s>\tau$ and $\e=\beta f'(s)(f'(s)-1)$. One has $\e>0$ and, coming back to \rf{a}, we obtain $f'''\leq-\e$ on $[s,+\infty)$.
After integration, we get 
$$\forall t\geq s,\qquad f''(t)-f''(s)\leq -\e (t-s)$$
and a contradiction with the fact that $f''>0$. 
Consequently, there exists $t_0\in(\tau,T_+)$ such that $f''>0$ on $[\tau,t_0)$ and $f''(t_0)=0$.
\ep

\bigskip

The last two lemmas give key results in the case where $\beta\in(0,1]$. The proofs can be found in \cite{rm} (see Lemma~5.16 and Lemma~A.11). However, for convenience, we give here proofs corresponding to the particular case that we consider.

\lem{\label{lf} Let $\beta\in(0,1]$ and $f$ be a solution of equation~\rf{eq} on some maximal interval $I=(T_-,T_+)$. If there exists $t_0\in I$ such that 
$$0< f'(t_0)<1\qquad\mbox{and}\qquad 0\leq f''(t_0)\leq f(t_0)(1-f'(t_0)),$$ then 
$T_+=+\infty$ and $f'(t)\to 1$ as $t\to+\infty$. Moreover, $f''>0$ on $[t_0,+\infty)$.}
\pr Let $\tau=\sup A(t_0)$ where 
$$A(t_0)=\left\{t\in[t_0,T_+)\;; \,f'(t_0)<f'<1\;\mbox{ and }\; f''>0 \;\mbox{ on }\; (t_0,t)\right\}.$$
The set $A(t_0)$ is not empty. This is clear if $f''(t_0)>0$, and if $f''(t_0)=0$ it follows from the fact that 
$f'''(t_0)=-\beta f'(t_0)(f'(t_0)-1)>0$. 
We claim that $\tau=T_+$.
Assume for contradiction that $\tau<T_+$. From Lemma~\ref{la}, item 2, we get that $f''(\tau)>0$, which implies, by definition of $\tau$, that $f'(\tau)=1$.
Therefore, since the function $H_f$ defined by~\rf{hf} is nonincreasing on $[t_0,\tau]$, we obtain
$$f''(\tau)=H_f(\tau)\leq H_f(t_0)=f''(t_0)+f(t_0)(f'(t_0)-1)\leq 0,$$
a contradiction. Thus, we have $\tau=T_+$. From Lemma~\ref{lb}, it follows that $T_+=+\infty$.
Since $f''>0$ on $[t_0,+\infty)$, by virtue of Lemma~\ref{ld}, we get that $f'(t)\to 1$ as $t\to+\infty$. \ep

\rem{\label{asym} If $f(t_0)>0$ and $f''(t_0)=0$, then $f(t)-t\to-\infty$ as $t\to+\infty$ (cf. \cite{rm}, Theorem 6.4, item 2.a).} 

\lem{\label{lff} Let $\beta\in(0,1]$ and $f$ be a solution of~\rf{eq} on some maximal interval $I=(T_-,T_+)$. If there exists $t_0\in I$ such that 
$$f'(t_0)>1\qquad\mbox{and}\qquad f(t_0)(1-f'(t_0))\leq f''(t_0)\leq 0,$$ 
then $T_+=+\infty$ and $f'(t)\to 1$ as $t\to+\infty$. Moreover, $f''<0$ on $[t_0,+\infty)$.}
\pr If we set $\tau=\sup B(t_0)$ where 
$$B(t_0)=\left\{t\in[t_0,T_+)\;; \,1<f'<f'(t_0)\;\mbox{ and }\; f''<0 \;\mbox{ on }\; (t_0,t)\right\},$$
the conclusion will follow by proceeding in the same way as in the previous proof.\ep

\rem{\label{asymp} If $f(t_0)>0$ and $f''(t_0)=0$, then $f(t)-t\to+\infty$ as $t\to+\infty$ (cf. \cite{rm}, Theorem 5.19, item 2.a).}

\section{\label{s3}Description of our approach when $b\geq 1$}

Let $\beta>0$, $a\geq0$ and $b\geq 1$. As said in the introduction, the method we will use to obtain solutions of the boundary value problems $({\cal P}_{\beta;a,b,0})$ and $({\cal P}_{\beta;a,b,1})$ is the shooting technique. 
Specifically, for $c\in\R$, let us denote by $f_c$ the solution of equation \rf{eq} satisfying the initial conditions 
\begin{equation}
f_c(0)= a, \quad\; f_c'(0)=b \quad \mbox
{ and }\quad f_c''(0)=c \label{ci}
\end{equation}
and let $[0,T_c)$ be the right maximal interval of existence of $f_c$. Hence, finding a solution of one of the problems 
$({\cal P}_{\beta;a,b,0})$ or $({\cal P}_{\beta;a,b,1})$ amounts to finding a value of $c$ such that $T_c=+\infty$ and 
$f'_c(t)\to 0$ or $1$ as $t\to+\infty$.
\medskip

To this end, let us partition $\R$ into the four sets ${\cal C}_0,\ldots,{\cal C}_3$ (or less if some of them are empty) defined as follows.
Let ${\cal C}_0=(0,+\infty)$ and, according to the notations used in \cite{rm}, let us set
\begin{align*}
&{\cal C}_1=\big\{c\leq 0\;;\;1\leq f'_c\leq b\;\mbox{ and }\;f''_c\leq 0\mbox{ on }[0,T_c)\big\} \\ \noalign{\vskip3mm}
&{\cal C}_2=\big\{c\leq 0\;;\;\exists\, t_c\in[0,T_c),\;\exists\, \epsilon_c>0  \mbox{ s.\,t. }f'_c>1\mbox{ on }(0,t_c),\\ \noalign{\vskip1,5mm}
&\hskip5cm f'_c<1 
\mbox{ on }(t_c,t_c+\epsilon_c)
\mbox{ and }f''_c<0\mbox{ on }(0,t_c+\epsilon_c)\big\}\\ \noalign{\vskip3mm}
&{\cal C}_3=\big\{c\leq 0\;;\;\exists\, r_c\in[0,T_c),\;\exists\, \eta_c>0 \mbox{ s.\,t. } f''_c<0 \mbox{ on }(0,r_c),\\ \noalign{\vskip1.5mm}
&\hskip5cm f''_c>0 \mbox{ on }(r_c,r_c+\eta_c)
\mbox{ and } f'_c>1\mbox{ on }(0,r_c+\eta_c) \big\}.
\end{align*}
This is obvious that ${\cal C}_0,\ldots,{\cal C}_3$ are disjoint sets and that their union is the whole line of real numbers. 

Thanks to Lemmas~\ref{lb} and \ref{ld}, if $c\in{\cal C}_1$ then $T_c=+\infty$ and $f'_c(t)\to1$ as $t\to+\infty$. 
In fact, ${\cal C}_1$ {\em is the set of values of $c$ for which $f_c$ is a concave solution of\,} $({\cal P}_{\beta;a,b,1})$.

\medskip

Since $\beta>0$, 
the study done in \cite{rm} (especially in Section 5.2) says, on the one hand, that ${\cal C}_3=\emptyset$ (which can easily be deduced from Lemma~\ref{la}, item~1) and, on the other hand, that either ${\cal C}_1=\emptyset$ and ${\cal C}_2=(-\infty,0]$, or there exists $c^*\leq 0$ such that ${\cal C}_1=[c^*,0]$ and ${\cal C}_2=(-\infty,c^*)$.  
In addition, if $\beta\in(0,1]$ then we are in the second case and $c^*\leq -a(b-1)$. If $\beta>1$ and $a=0$  then
${\cal C}_1=\emptyset$, but, for $a>0$, we do not know if ${\cal C}_1$ is empty or not. 

\medskip

In the next sections we will distinguish between the cases $\beta\in(0,1]$ and $\beta>1$. In the first case, we can give a complete description of the solutions (see Theorem~\ref{th}), whereas in the second one, we have only partial answers.

We will also consider the case where $b\in[0,1)$, for which we will have to partition $\R$ in a slightly different way.  

\medskip

Before that, and in order to complete the study, let us divide the set ${\cal C}_2$ into the following two subsets 
\begin{align*}
&{\cal C}_{2,1}=\{c\in{\cal C}_2\;;\;f'_c>0\mbox{ on } [0,T_c) \}\\ \noalign{\vskip2mm}
&{\cal C}_{2,2}=\{c\in{\cal C}_2\;;\;\exists\, s_c\in (0,T_c)\;\mbox{ s.\,t. } f'_c>0 \mbox{ on }[0,s_c)
\mbox{ and } f'_c(s_c)=0\}
\end{align*}
and let us give properties of each of them that hold for all $\beta>0$.

\lem{\label{pj} If $c\in\R$ is such that $f'_c>0$ on $[0,T_c)$, then $T_c=+\infty$ and $\,f'_c$ is bounded.
Moreover, if $c\leq0$, then $\,f'_c\leq\max\{b\,;\frac32\}$ on $[0,+\infty)$.
}
\pr 
Let $c\in\R$ be such that $f'_c>0$ on $[0,T_c)$. Then $f_c\geq a\geq0$ on $[0,T_c)$, and thanks to Lemma~\ref{lee},
it follows that $T_c=+\infty$ and that $f'_c$ is bounded.

It remains to show that $f'_c\leq\max\{b\,;\frac32\}$ in the case where $c\leq0$. As in \rf{fonctionlf}, let us define the function $L_c$  on $[0,+\infty)$ by
\begin{equation}
L_c(t)=3f''_{c}(t)^2+\beta (2f'_{c}(t)-3)f'_{c}(t)^2.\label{fonctionl}
\end{equation}
We have $L_c'(t)=-6f_c(t)f_c''(t)^2$ and, since $f_c\geq0$, it implies that $L_c$ is nonincreasing.

If $f''_c\leq0$ on $(0,+\infty)$, then $f'_c\leq b$. Otherwise, there exists $t_0$ such that $f''_c<0$ on $(0,t_0)$ and
$f''_c(t_0)=0$ (which can occur only when $c<0$, or $c=0$ and $b>1$). By Lemma~\ref{la}, item~3, it follows that $f'_c(t_0)<1$, and thus $L_c(t_0)<0$. Then, $L_c<0$ on $(t_0,+\infty)$ which implies that $f'_c\leq \frac32$ on $(t_0,+\infty)$. Since $f'_c\leq b$ on $(0,t_0)$, the proof is complete.
\ep

\propo{\label{pk} Let $c_*=\sup\,({\cal C}_{1}\cup{\cal C}_{2,1})$. Then $c_*$ is finite.
}
\pr  
Let $c\in{\cal C}_{1}\cup{\cal C}_{2,1}$. By the definition of ${\cal C}_{1}$ and ${\cal C}_{2,1}$, and thanks to Lemma~\ref{pj}, we have $T_c=+\infty$ and $0<f'_c\leq d$  on $(0,+\infty)$ where $d=\max\{b\,;\frac32\}$.

Since $(f''_c+f_cf'_c)'=-\beta f'_c(f'_c-1)+{f'_c}^2\leq\beta f'_c+{f'_c}^2\leq d(\beta +d)$, by integrating, we then have
$$
\forall t\geq0,\qquad f''_c(t)+f_c(t)f'_c(t)\leq c+ab+d(\beta +d) t.
$$
Integrating once again, we get
$$
\forall t\geq0,\qquad 0<f'_c(t)\leq f'_c(t)+\tfrac12f_c(t)^2\leq b+\tfrac12{a^2}+(c+ab)t+\tfrac12d(\beta +d) t^2$$
which implies that $c\geq-ab-\sqrt{(2b+a^2)(\beta+d)d}$.
\ep

\rem{As we have seen above, if ${\cal C}_{1}\neq\emptyset$, then ${\cal C}_{1}=[c^*,0]$ and thus ${\cal C}_{2,1}\subset [c_*,c^*)$. }

\propo{\label{pl} We have $(-\infty,c_*)\subset{\cal C}_{2,2}$. Moreover, if $c\in{\cal C}_{2,2}$ then $T_c<+\infty$ and $f''_c<0$ on $(0,T_c)$.
}
\pr The fact that $(-\infty,c_*)\subset{\cal C}_{2,2}$ follows immediately from Proposition~\ref{pk}. 
Let $c\in{\cal C}_{2,2}$. 
Then, there exists $s_c\in (0,T_c)$ such that $f'_c>0$ on $[0,s_c)$ and $f'_c(s_c)=0$. Consider the function $L_c$ defined by \rf{fonctionl}. Since $f_c\geq 0$ on $[0,s_c]$, then $L_c$ is nonincreasing on $[0,s_c]$. 

Suppose first that $c<0$. Assume for contradiction that there exists $t_0\in[0,s_c)$ such that 
$f''_c<0$ on $[0,t_0)$ and $f''_c(t_0)=0$, then $0<f'_c(t_0)<1$ (see Lemma~\ref{la}, item 3), and hence $L_c(t_0)<0$. Since $L_c$ is nonincreasing on $[0,s_c]$, this contradicts the fact that 
$L_c(s_c)=3f''_c(s_c)^2\geq 0$. Therefore, $f''_c<0$ on $[0,s_c]$. 

If $c=0$, which can only happen if $b>1$, then $f'''_c(0)=-\beta b(b-1)<0$. Hence there exists $\eta\in(0,s_c)$ such that $f''_c<0$ and $f'_c>1$ on $(0,\eta]$. The arguments above applied to the function $t\mapsto f_c(t+\eta)$ give that $f''_c<0$ on $[\eta,s_c]$ and thus on $(0,s_c]$.

To get that $f''_c<0$ on $(0,T_c)$, it remains to notice that $f''_c$ cannot vanish on $(s_c,T_c)$, by virtue of Lemma~\ref{la}, item~3.

Finally, the fact that $T_c<+\infty$ follows from Proposition~2.11 of \cite{rm}, which says that, for any $\tau\in\R$, there is no negative (strictly) concave function $f$ such that $f'''+ff''\leq 0$ on $[\tau,+\infty)$.
\ep

\rem{\label{rl} If $c\in{\cal C}_{2,2}$ then $f_c$ is strictly concave on $[0,T_c)$, has a global maximum at $s_c$ and $f_c(t)\to-\infty$ as $t\to T_c$. In addition, $f'_c(t)$ and $f''_c(t)$ tend to $-\infty$ as $t\to T_c$.}

\propo{\label{pm} The set ${\cal C}_{2,2}$ is an open set of $(-\infty,0]$ $($for its induced topology$)$.}
\pr 
Let $c_{0}\in {\cal C}_{2,2}$. There exists $\tau \in (0,T_{c_{0}})$ such that 
$f'_{c_{0}}(\tau )<0$. Let us set $\varepsilon =-\frac12 f'_{c_{0}}(\tau )$.
By continuity of the function $c\mapsto f'_c(\tau)$, there exists
$\alpha >0$ such that, for all $c\in (-\infty,0]$, one has
$$\vert c-c_0\vert <\alpha \;\Longrightarrow\; f'_c(\tau )<f'_{c_{0}}(\tau )+\varepsilon.$$
Therefore, $f'_{c}(\tau )<0$ and $c\in{\cal C}_{2,2}$.
\ep

\section{\label{s4}The case $\beta\in(0,1]$ and $b\geq 1$}

In this section, we assume that $\beta\in(0,1]$, $a\geq0$ and $b\geq 1$.

\propo{\label{pn} If $c\in{\cal C}_{0}$, then $T_c=+\infty$ and $f'_c(t)\to1$ as $t\to+\infty$.
}
\pr From Lemma~\ref{pi}, there exists $t_0\in(0,T_c)$ such that $f_c''>0$ on $[0,t_0)$ and $f_c''(t_0)=0$. Since $f_c(t_0)>0$ and $f_c'(t_0)> b>1$, the conclusion follows from Lemma~\ref{lff}.
\ep

\rem{Thanks to the previous proposition, we see that $f_c$ is a {\em convex-concave} of $({\cal P}_{\beta;a,b,1})$ for all $c>0$.
Moreover, we have that
$f_c(t)-t\to+\infty$ as $t\to+\infty$ (cf. Remark~\ref{asymp}). }

\propo{\label{pz} There exists $c^*\leq -a(b-1)$ such that ${\cal C}_{1}=[c^*,0]$.}
\pr If $b=1$ then ${\cal C}_{1}=\{0\}$. If $b>1$, as we already said in the previous section, this result is proven in \cite{rm} (see Corollary~5.13 and Lemma~5.16). For convenience, let us recall briefly the main arguments which were used to get it.
On the one hand, from Lemma~\ref{lff} with $t_0=0$ (or Lemma~5.16 of \cite{rm}), it follows that $[-a(b-1),0]\subset{\cal C}_{1}$. On the other hand, Lemma~5.12 of \cite{rm} implies that ${\cal C}_2$ is an interval of the type $(-\infty,c^*)$. This completes the proof since ${\cal C}_1=(-\infty,0]\setminus{\cal C}_2$.
\ep 

\rem{\label{rko} From the previous proposition, we have that $0\notin {\cal C}_{2,2}$. Hence, Proposition~\ref{pm} implies that ${\cal C}_{2,2}$ is an open set.}

\propo{\label{pa} If $c\in {\cal C}_{2,1}$ then $T_c=+\infty$ and $f'_{c}$ has a finite limit 
at infinity, equal either to $0$ or to $1.$
}
\pr Let $c\in {\cal C}_{2,1}$. By Proposition~\ref{pz}, we have $c<0$.
Thanks to Lemma~\ref{pj}, we know that $T_c=+\infty$.
Assume first that $f''_c<0$ on $(0,+\infty)$. 
Then $f'_c$ is positive and decreasing, and thus $f'_c$ has a finite limit $\lambda\geq0$ at infinity. Moreover, $f'_c$ takes the value $1$ at some point, hence $\lambda\in[0,1)$ and, by Lemma~\ref{ld}, we finally get that $\lambda=0$.

Assume now that $f''_c$ vanishes on $(0,+\infty)$. Let $t_0$ be the first point where $f''_c$ vanishes. Thanks to 
Lemma~\ref{la}, item~3, we have  
$0<f'_c(t_0)<1$, and the conclusion follows from Lemma~\ref{lf}.
\ep

\rem{\label{rcc} If $c\in {\cal C}_{2,1}$ then either $f_c$ is a {\em concave\,} solution of $({\cal P}_{\beta;a,b,0})$ or $f_c$ is a {\em concave-convex\,} solution of $({\cal P}_{\beta;a,b,1})$. In the first case, there exists $l>a$ such that $f_c(t)\to l$ as $t\to+\infty$ (cf. Lemma~\ref{le}) and, in the second one, we have that $f_c(t)-t\to-\infty$ as $t\to+\infty$ (cf. Remark~\ref{asym}).}

\propo{\label{pg}
Let $c\in {\cal C}_{2,2}$. For all $t\in[0,T_c)$, one has $f_{c}(t)\leq \sqrt{a^2+2b}$.
}
\pr
Let $c\in {\cal C}_{2,2}$ and $s_c$ be as in the definition of ${\cal C}_{2,2}$, i.e. such that $f'_c>0$ on $[0,s_c)$ and $f'_c(s_c)=0$. 
For all $t\in[0,s_c]$, we have 
\begin{align}
\nonumber\big(tf''_c(t)-f'_c(t)+tf_c(t)f'_c(t)\big)'
&= tf'''_c(t)+tf_c(t)f''_c(t)+tf'_c(t)^2+f_c(t)f'_c(t)\\ \noalign{\vskip2mm}
&= (1-\beta)tf'_c(t)^2+\beta tf'_c(t)+f_c(t)f'_c(t)\geq f_c(t)f'_c(t).\label{ineqb}
\end{align}
Integrating between $0$ and $s_c$ yields
$$f_c(s_c)^2\leq a^2+2\big(s_cf''_c(s_c)+b\big)\leq a^2+2b$$
and
$f_{c}(s_c)\leq \sqrt{a^2+2b}$.
The conclusion follows from the fact that,
for all $t\in[0,T_c)$, we have $f_c(t)\leq f_c(s_c)$, as we noticed in Remark~\ref{rl}.
\ep

\propo{\label{ph}
Let $c$ be a point of the boundary of \,${\cal C}_{2,2}$. Then, $c\in{\cal C}_{2,1}$ and $f'_c(t)\to0$ as $t\to+\infty$.
Moreover, $f_{c}$ is bounded and concave.
}
\pr
Let $c$ be a point of the boundary of ${\cal C}_{2,2}$ and $(c_n)_{n\geq0}$ be a sequence of ${\cal C}_{2,2}$ such that $c_n\to c$ as $n\to+\infty$. For all $n\geq 0$, let us set $T_n=T_{c_n}$ and $f_n=f_{c_n}$.
Since ${\cal C}_{2,2}$ is an open set, then $c\in{\cal C}_{1}\cup {\cal C}_{2,1}$ and hence $T_c=+\infty$.
Let $t\geq 0$ be fixed. From the lower semicontinuity of the function $d\to T_d$, we get that there exists $n_0\geq0$ such that $T_n\geq t$ for all $n\geq n_0$. Since $f_n(t)\to f_c(t)$ as $n\to+\infty$, we deduce from Proposition~\ref{pg} that 
$f_c$ is bounded. Therefore, $f_c'$ cannot tend to $1$ at infinity and thus, necessarily, we have $c\in{\cal C}_{2,1}$ and $f'_c(t)\to0$ as $t\to+\infty$. Moreover, $f_c$ is concave (cf. Remark~\ref{rcc}).
\ep

\propo{\label{pu} There exists at most one $c$ such that $f'_c(t)\to0$ as $t\to+\infty$.
}
\pr
From Proposition~\ref{pz}, Proposition~\ref{pa} and Lemma~\ref{le}, we see that if $c$ is such that $f'_c(t)\to0$ as $t\to+\infty$, then $c<0$,
$f''_c<0$ and $f_c$ is bounded. 
For such a $c$, as done in \cite{rm},~Section~4, we can define a function $v:(0,b^2]\rightarrow\R$ such that
\begin{equation}
\forall t\geq0,\quad v(f_c'(t)^2)=f_c(t).\label{vy}
\end{equation}
By setting $y=f_c'(t)^2$, we get 
\begin{equation*}
f_c(t)=v(y),\quad f_c'(t)=\sqrt{y},\quad f''_c(t)=\frac{1}{2v'(y)}\quad \text{and}\quad f_c'''(t)=-\frac{v''(y)\sqrt{y}}{2v'(y)^3}\label{y1}
\end{equation*}
and using \rf{eq} we obtain
\begin{equation}
\forall y\in (0,b^2],\quad v''(y)=\frac{v(y)\,v'(y)^2}{\sqrt{y}}+2\beta(\sqrt{y}-1)\,v'(y)^3. \label{eqv}
\end{equation}  
From \rf{ci}, we deduce that $v(b^2)=a$ and $v'(b^2)=\frac1{2c}$. Moreover, since $f_c$ is bounded, it is so for $v$.

Assume that there exists $c_1>c_2$ such that $f'_{c_1}(t)\to0$ and $f'_{c_2}(t)\to0$ as $t\to+\infty$, and denote by $v_1$ and $v_2$ the functions associated to $f_{c_1}$ and $f_{c_2}$ by \rf{vy}.  If we set $w=v_1-v_2$ then $w(b^2)=0$ and $w'(b^2)<0$. We claim that $w'<0$ on $(0,b^2]$. 
For contradiction, assume there exists $x\in(0,b^2)$ such that $w'<0$ on $(0,x)$ and $w'(x)=0$. Hence we have $w''(x)\leq 0$ and $w(x)>0$. But, thanks to \rf{vy}, we have
$$w''(x)=\frac{w(x)}{\sqrt x}\,v'_1(x)^2$$
and a contradiction. 

Now, let us set $V_i=1/{v'_i}$ for $i=1,2$ and $W=V_1-V_2$. Then $W(b^2)=2(c_1-c_2)>0$ and $W(y)\to0$ as $y\to 0$. In the other hand, thanks to \rf{eqv}, we have
$$\forall y\in(0,b^2],\qquad W'(y)=-\frac{w(y)}{\sqrt y}-2\beta(\sqrt{y}-1)\,w'(y).$$
Therefore,  we have 
\begin{align}
\nonumber W(b^2)&=\int_0^{b^2} W'(y)\,dy=-\int_0^{b^2} \bigg(\frac{w(y)}{\sqrt y}+2\beta(\sqrt{y}-1)\,w'(y)\bigg)dy\\ \noalign{\vskip2mm}
\nonumber &=-2\Big[\sqrt y\,w(y)\Big]_0^{b^2}+2\int_0^{b^2}\big((1-\beta)\sqrt y+\beta\big)\,w'(y)\,dy\\ \noalign{\vskip2mm}
&=2\int_0^{b^2}\big((1-\beta)\sqrt y+\beta\big)\,w'(y)\,dy,\label{wv}
\end{align}
the last equality following from the fact that $w(y)$ tends to a finite limit as $y\to0$.
Since $w'<0$, we finally obtain $W(b^2)<0$ and a contradiction.
\ep

\rem{\label{rv} The change of variable \rf{vy} is particularly efficient to obtain some uniqueness results. 
In \cite{rm}, it is used for the general equation $f'''+ff''+\g(f')=0$ (cf. Section 4, Lemma~5.4 and Lemma~5.17). The case we examined in Proposition~\ref{pu} is part of Lemma~5.17 of \cite{rm} with $\lambda=0$.
In this lemma, it is assumed that $0<\g(x)\leq x^2$ for $x\in(0,b]$ to ensure uniqueness.
Here, in Proposition~\ref{pu}, we have $\g(x)=\beta x(x-1)$ with $\beta\in(0,1]$ and hence $\beta x(x-1)\leq x^2$ for $x\in(0,b]$, but $\beta x(x-1)\leq 0$ for $x\in(0,1]$. However, the assumption about the positivity of $\g$ is not relevant 
because not used in the proof of Lemma~5.17 of \cite{rm}. 
In addition, the inequality $\beta x(x-1)\leq x^2$ is still true on $(0,b]$, if $\beta>1$ and $1\leq b\leq\frac\beta{\beta-1}$. Finally, let us notice that, in the latter case, the integral in~\rf{wv} is still negative, and the contradiction occurs there too.
}

\coro{\label{corcc} One has  $\,{\cal C}_{2,2}=(-\infty,c_*)$ and $\,{\cal C}_{2,1}=[c_*,c^*)$.}
\pr From Remark~\ref{rko}, Propositions~\ref{pl},~\ref{ph} and~\ref{pu}, we see that ${\cal C}_{2,2}$ is open, contains $(-\infty,c_*)$ and its boundary is reduced to a single point. Therefore, since $c_*=\sup\,({\cal C}_{1}\cup{\cal C}_{2,1})$, we necessarily have ${\cal C}_{2,2}=(-\infty,c_*)$ and ${\cal C}_{2,1}=[c_*,c^*)$.
\ep

To finish this section, let us express the results of Proposition~\ref{pn}, Proposition~\ref{pz} and Corollary~\ref{corcc} in terms of the boundary problems $({\cal P}_{\beta;a,b,0})$ and $({\cal P}_{\beta;a,b,1})$.

\theo{\label{th} 
Let $\beta\in(0,1]$, $a\geq 0$ and $b\geq 1$. There exists $c_*<0$ such that $:$
\vspace{-1mm}
\bi  \itemsep=0mm
\item[$\triangleright$] $f_{c}$ is not defined on the whole interval $[0,+\infty)$ if $c<c_*$ $;$
\item[$\triangleright$] $f_{c_*}$ is a {\em concave} solution of $({\cal P}_{\beta;a,b,0})$ $;$
\item[$\triangleright$] $f_{c}$ is a solution of $({\cal P}_{\beta;a,b,1})$ for all $c\in(c_*,+\infty)$.
\ei
Moreover, there exists $c^*\in (c_*,-a(b-1)]$ such that $:$
\vspace{-1mm}
\bi  \itemsep=0mm
\item[$\triangleright$] $f_{c}$ is a {\em convex-concave} solution of $({\cal P}_{\beta;a,b,1})$ for all $c\in(0,+\infty)$ $;$
\item[$\triangleright$] $f_{c}$ is a {\em concave} solution of $({\cal P}_{\beta;a,b,1})$ for all $c\in[c^*,0]$ $;$
\item[$\triangleright$] $f_{c}$ is a {\em concave-convex} solution of $({\cal P}_{\beta;a,b,1})$ for all $c\in(c_*,c^*)$.
\ei
}

\medskip

\rem{The previous theorem says that problem $({\cal P}_{\beta;a,b,0})$ has one and only one solution, whereas problem $({\cal P}_{\beta;a,b,1})$ has infinite number of solutions.}

\rem{We know that $f_{c_*}$ has a finite limit at infinity, denoted by $l$. By slightly modifying the proof of Proposition~7.2 of \cite{rm}, one can prove that there exists a positive constant $A$ such that, for all $\epsilon>0$, the following hold
$$f_{c_*}''(t)=-l^2Ae^{-l t}\left(1+{\rm o}\!\left(e^{-(l-\epsilon)t}\right)\right),\qquad
f_{c_*}'(t)=l Ae^{-l t}\left(1+{\rm o}\!\left(e^{-(l-\epsilon)t}\right)\right)$$
$$f_{c_*}(t)=l-Ae^{-l t}\left(1+{\rm o}\!\left(e^{-(l-\epsilon)t}\right)\right)\quad\mbox{ as }\quad t\to+\infty.$$
}

\rem{Among the concave solutions of $({\cal P}_{\beta;a,b,1})$, only $f_{c^*}$ has a slant asymptote, i.e. there exists $l>a$ such that $f_{c^*}(t)-t\to l$ as $t\to+\infty$. In addition, Proposition~7.5 of \cite{rm} implies that, as $t\to+\infty$, we have
$$f_{c^*}''(t)=-e^{-\frac{t^2}{2}-l t+{\rm O}\left(\ln t\right)},\quad f_{c^*}'(t)=1+e^{-\frac{t^2}{2}-l t+{\rm O}\left(\ln t\right)}\quad\mbox{and}\quad f_{c^*}(t)= t+l-e^{-\frac{t^2}{2}-l t+{\rm O}\left(\ln t\right)}.$$
If $c^*<0$, then the function $t\mapsto f_c(t)-t$ is unbounded, for any $c\in(c^*,0]$. 

It is possible to do better and to precise what is the term ${\rm O}(\ln t)$. By a method used for the Falkner-Skan equation in \cite{hart}, Chapter XIV, Theorem 9.1, one can show that there exists a constant $A>0$ such that
$$f_{c^*}'(t)-1\sim At^{\beta-1}e^{-\frac{t^2}{2}-l t}\quad\mbox{ as }\quad t\to+\infty.$$ 
Other asymptotic results for $f_c$ (concave, convex-concave or concave-convex) such that the function $t\mapsto f_c(t)-t$ is unbounded, should also be obtained by applying the ideas of \cite{hart}, Chapter XIV, Theorem 9.1 and 9.2. See also \cite{sc}.
}

\rem{\label{outl} The main ingredients used in this section are, one the one hand, Lemmas~\ref{lf} and~\ref{lff} that precise the behavior of $f_c$ after a point where $f''_c$ vanishes and, on the other hand, the fact that the set ${\cal C}_{2,2}$ has at most one point on its boundary, implying that it is an interval.
}

\section{\label{s5}The case $\beta\in(0,1]$ and $0\leq b<1$}

Let $\beta\in(0,1]$, $a\geq0$ and $0<b<1$. In this situation, it is easy to see that $\R$
can be partitioned into the four sets ${\cal C}'_{0,1}$, ${\cal C}'_{0,2}$, ${\cal C}'_1$ and ${\cal C}'_2$ where
\begin{align*}
&{\cal C}'_{0,1}=\{c<0\;;\;f'_c>0\mbox{ on } [0,T_c) \}\\ \noalign{\vskip2mm}
&{\cal C}'_{0,2}=\{c<0\;;\;\exists\, s_c\in (0,T_c)\;\mbox{ s.\,t. } f'_c>0 \mbox{ on }[0,s_c)
\mbox{ and } f'_c(s_c)=0\}\\ \noalign{\vskip2mm}
&{\cal C}'_1=\big\{c\geq 0\;;\;b\leq f'_c\leq 1\;\mbox{ and }\;f''_c\geq 0\mbox{ on }[0,T_c)\big\} \\ \noalign{\vskip3mm}
&{\cal C}'_2=\big\{c\geq 0\;;\;\exists\, t_c\in[0,T_c),\;\exists\, \epsilon_c>0  \mbox{ s.\,t. }f'_c<1\mbox{ on }(0,t_c),\\ \noalign{\vskip1,5mm}
&\hskip5cm f'_c>1 
\mbox{ on }(t_c,t_c+\epsilon_c)
\mbox{ and }f''_c>0\mbox{ on }(0,t_c+\epsilon_c)\big\}.
\end{align*}
The fact that any $c\geq0$ belongs to ${\cal C}'_1\cup {\cal C}'_2$ is due inter alia to Lemma~\ref{la}, item~4, which implies that $f''_c$ remains positive as long as $f'_c\leq1$.

\smallskip
The arguments used in the previous section, and evoked in Remark~\ref{outl}, can be applied here. Some results, as Propositions~\ref{pg} and~\ref{ph}, are still true. On the other hand, as we will see below, some other results are obtained more easily. For example, the existence and the uniqueness of a concave solution of $({\cal P}_{\beta;a,b,0})$ are already known, and so it is not necessary to argue as in the previous section (cf. Propositions~\ref{ph} and~\ref{pu}).
\medskip

Since $\beta x(x-1)<0$ for $x\in(0,b]$, it follows from Theorem 5.5 of \cite{rm} that there exists a unique $c_*$ such that $f_{c_*}$ is a concave solution of $({\cal P}_{\beta;a,b,0})$. Moreover, we have $c_*<0$. As in the previous section, this implies that ${\cal C}'_{0,2}=(-\infty,c_*)$. Hence ${\cal C}'_{0,1}=[c_*,0)$, and if $c\in(c_*,0)$, then $f''_c$ vanishes at a first point where $f'_c<1$.

Next, in the same way as in the proof of Proposition~\ref{pk}, we can prove that $c^*=\sup {\cal C}'_1$ is finite, and hence that ${\cal C}'_1=[0,c^*]$ and ${\cal C}'_2=(c^*,+\infty)$. 
Moreover, from Lemma~\ref{lf}, we have $c^*\geq a(1-b)$. 
On the other hand, it follows from Lemma~\ref{pi} that, if $c\in{\cal C}'_2$, then $f''_c$ vanishes at a first point where $f'_c>1$.

All this, combined with an appropriate use of Lemmas~\ref{lf} and~\ref{lff}, allows to state the following theorem. For more details, we refer to \cite{bk}.
\medskip

\theo{\label{thh} Let $\beta\in(0,1]$, $a\geq 0$ and $b\in(0,1)$. There exist $c_*<0$ and $c^*\geq a(1-b)$ such that $:$
\vspace{-2mm}
\bi  \itemsep=0mm
\item[$\triangleright$] $f_{c}$ is not defined on the whole interval $[0,+\infty)$ if $c<c_*$ $;$
\item[$\triangleright$] $f_{c_*}$ is a {\em concave} solution of $({\cal P}_{\beta;a,b,0})$ $;$
\item[$\triangleright$] $f_{c}$ is a {\em concave-convex} solution of $({\cal P}_{\beta;a,b,1})$ for all $c\in(c_*,0)$ $;$
\item[$\triangleright$] $f_{c}$ is a {\em convex} solution of $({\cal P}_{\beta;a,b,1})$ for all $c\in[0,c^*]$ $;$
\item[$\triangleright$] $f_{c}$ is a {\em convex-concave} solution of $({\cal P}_{\beta;a,b,1})$ for all $c\in(c^*,+\infty)$.
\ei
}

\smallskip

\rem{\label{bo} [\,{\sc The case $b=0$}\,] We can show similar results if $b=0$. For details of the proof, we refer to \cite{bk}. 
\vspace{-2mm}
\bi  \itemsep=0mm
\item[$\triangleright$] If $c<0$, then $T_c<+\infty$. 
\item[$\triangleright$] For $c=0$, we have $f_0(t)=a$. 
\item[$\triangleright$] There exists $c^*\geq a$ such that $f_{c}$ is a convex solution of $({\cal P}_{\beta;a,b,1})$ for all $c\in(0,c^*]$ and $f_{c}$ is a convex-concave solution of $({\cal P}_{\beta;a,b,1})$ for all $c>c^*$.\ei
}

\section{\label{s6}About the case $\beta>1$}

In this section, we will assume that $\beta>1$, $a\geq0$ and $b>0$. The main difference with the case $\beta\in(0,1]$, is that Lemmas~\ref{lf} and~\ref{lff} do not necessarily hold anymore. In fact, it is the case if $f(t_0)=0$, and in particular this implies that, if $a=0$ and $b>1$, then we have ${\cal C}_1=\emptyset$ (see \cite{rm}, Theorem~5.19, item~2.b), and  if $a=0$ and $0<b<1$, then ${\cal C}'_1=\emptyset$
(see \cite{rm}, Theorem~6.4, item~2.b).

Another consequence is that, on the contrary to what happens in the case $\beta\in(0,1]$, where for any $c$ the function $f''_c$ vanishes at most once in $[0,T_c)$, this is not necessarily true if $\beta>1$, and numerical experimentations indicate that it is so.

Furthermore, nothing indicates whether 
both problems $({\cal P}_{\beta;a,b,0})$ and $({\cal P}_{\beta;a,b,1})$ have solutions or not.

Nevertheless, some results are still true.
We start with a result about the problem $({\cal P}_{\beta;a,b,0})$. Next, we prove that, if $f'_c$ remains positive, then $f'_{c}$ tends to $0$ or $1$ at infinity.
Finally, we point some situations for which the problem $({\cal P}_{\beta;a,b,1})$ has solutions.

\propo{\label{paaa} If $b\in(0,\frac{\beta}{\beta-1}]$, then there exists $c_*<0$ such that $f_{c_*}$ is a solution of the problem $({\cal P}_{\beta;a,b,0})$. 
Moreover, $f_{c_*}$ is concave and is the unique solution of $({\cal P}_{\beta;a,b,0})$.
}
\pr If $b\in(0,1)$, as in the previous section, this follows from \cite{rm}, Theorem 5.5. If $b\in[1,\frac{\beta}{\beta-1}]$, on the one hand, we remark that inequality~\rf{ineqb} still holds, and hence it is so for the conclusions of Propositions~\ref{pg} and~\ref{ph}. Thus, the problem $({\cal P}_{\beta;a,b,0})$ has a solution. On  the other hand, as we point out in Remark~\ref{rv}, the uniqueness of the solution of $({\cal P}_{\beta;a,b,0})$ holds true for $b\in[1,\frac{\beta}{\beta-1}]$. 
\ep

\propo{\label{paa} If $c\in \R$ is such that $f'_c>0$ on $(0,T_c)$, then $T_c=+\infty$ and $f'_{c}$ has a finite limit 
at infinity, equal either to $0$ or to $1$.
}
\pr Let $c\in\R$ be such that $f'_c>0$ on $(0,T_c)$.
From Lemma~\ref{lee}, we know that $T_c=+\infty$ and that $f'_c$ is bounded. 

If there exists a point $\tau\geq 0$ such that $f''_c$ does not change of sign on $(\tau,+\infty)$, then $f'_c$ is monotone on this interval. Hence, $f'_c$ has a finite limit at infinity and, by virtue of Lemma~\ref{ld}, this limit is equal to $0$ or $1$.

If we are not in the previous situation, then there exists an increasing sequence $(\tau_n)_{n\geq0}$ tending to $+\infty$ such that $f''_c(\tau_n)=0$ and $f'''_c(\tau_n)>0$, for all $n\geq0$ (notice that Lemma~\ref{lc} implies that we cannot have $f'''_c(\tau_n)=0$). 

Let $L_c$ be the function defined on $[0,+\infty)$ by~\rf{fonctionl}, i.e. 
$$\forall t\geq 0,\qquad L_c(t)=3f''_{c}(t)^2+\beta (2f'_{c}(t)-3)f'_{c}(t)^2.$$
We know that $L_c$ is decreasing and takes negative value at each $\tau_n$ since, by virtue of Lemma~\ref{la}, item 3, we have $f'_c(\tau_n)<1$.
Therefore, we have  $L_c(t)<0$ for $t\geq \tau_0$. Moreover, since $2x^3-3x^2\geq -1$ for $x\geq0$, then $L_c(t)\geq-\beta$ for all $t\geq0$. Hence $L_c(t)$ tends to some $\alpha<0$ as $t\to+\infty$.

Inspired by an idea developed in \cite{guedda} we will show that $f_c(t)\to +\infty$ and $f''_c(t)\to0$ as $t\to+\infty$.

First, let us prove that $f_c(t)\to+\infty$ as $t\to+\infty$. If it is not the case, then $f_c$ has a finite limit $l$ at infinity (recall that $f_c$ is increasing) and there exists a sequence $(s_n)_{n\geq0}$ in $[\tau_0,+\infty)$ such that $s_n\to+\infty$ and $f'_c(s_n)\to 0$ as $n\to +\infty$. 

By passing to the limit as $n\to+\infty$ in the inequalities 
$$\beta f'_c(s_n)^2(2f'_c(s_n)-3)\leq L_c(s_n)\leq L_c(\tau_0)<0$$
we get a contradiction. Therefore $f_c(t)\to+\infty$ as $t\to+\infty$.

Next, let us prove that $f''_c(t)\to0$ as $t\to+\infty$.  
Let $x_n$ be a point of the interval $(\tau_n,\tau_{n+1})$ such that $\vert f''_c(t)\vert\leq \vert f''_c(x_n)\vert$ for all $t\in[\tau_n,\tau_{n+1}]$. We have $f'''_c(x_n)=0$ and thus,
from equation~\rf{eq}, one has 
\begin{equation*}
f''_{c}(x_{n})=\frac{-\beta f'_{c}(x_{n})(f'_{c}(x_{n})-1)}{f_{c}(x_{n})}.
\end{equation*}
Thus, since $f'_c$ is bounded and that $f_c(x_n)\to+\infty$ as $n\to+\infty$, we obtain that $f''_c(x_n)\to0$ as $n\to+\infty$, and hence $f''_c(t)\to0$ as $t\to+\infty$.

Now we are able to conclude. 
Since $f''_{c}(t)\to 0$ and $L_c(t)\to\alpha$ as $t\to+\infty$, we have 
that $2f'^3_{c}(t)-3f'^2_{c}(t)\to\alpha$ as $t\to+\infty$.
Therefore $f'_c$ has a finite limit $\lambda$ at infinity, that is a root of the polynomial $P(x)=2x^3-3x^2-\alpha$ (see Remark~\ref{raa} below).
Since $P(0)=-\alpha\neq0$, by Lemma~\ref{ld}, we get $\lambda=1$.
\ep

\rem{\label{raa} In the previous proof, we used the fact that for any real polynomial $P$ with real roots $a_1,\ldots,a_s$ and any continuous function $\f:[0,+\infty)\to\R$ such that $P(\f(t))\to0$ as $t\to+\infty$, then $\f(t)$ tends to a root of  $P$ as 
$t\to+\infty$.
To prove this, note first that, for every $\e$ small enough, the intervals $A_{j,\e}=]a_j-\e,a_j+\e[$ are disjoint. Denote by  $A_\e$ their union.
On the one hand, since $P(\f(t))\to0$ as $t\to+\infty$, for all $n\geq1$, there exists $t_n$ such that $\f([t_n,+\infty[\,)\subset P^{-1}([-\frac1n,\frac1n])$.
On the other hand, since 
$$\bigcap_{n\geq1}P^{-1}\big(\big[-\tfrac1n,\tfrac1n\big]\big)
=P^{-1}(\{0\})=\{a_1,\ldots,a_s\},$$
there exists $n_\e$ such that $P^{-1}([-\frac1{n_\e},\frac1{n_\e}])\subset A_\e$. Set $t_\e=t_{n_\e}$, one has 
$\f([t_\e,+\infty[\,)\subset A_\e$. Due to the continuity of $\f$ the set $\f([t_\e,+\infty[\,)$ is an interval, and hence there exists $k\in\{1,\ldots,s\}$ such that $\f([t_\e,+\infty[\,)\subset A_{k,\e}$. In other words,
for $t\geq t_\e$ we have $\vert \f(t)-a_k\vert<\e.$
Finally, $\f(t)\to a_k$ as 
$t\to+\infty$.
}

\rem{\label{rab} In the proof of Proposition~\ref{paa}, we only use the positivity of $\beta$. Thus Proposition~\ref{paa}
implies Proposition~\ref{pa}, but the proof of this latter proposition is simpler and shorter, and says more, i.e. that $f''_c$ vanishes at most once.
}

\propo{\label{pab} If $\beta\in(1,2]$ and $a>0$, then for any $c$ such that $2ac\geq b^2-(2b-\beta)a^2$, we have $T_c=+\infty$ and $f'_c(t)\to 1$ as $t\to+\infty$.
}
\pr Let $c\in\R$ and denote by $K_c$ the function defined on $[0,T_c)$ by 
\begin{equation*}
K_c(t)=2f_c(t)f''_c(t)-f'_c(t)^2+(2f'_c(t)-\beta)f_c(t)^2.
\end{equation*}
From \rf{eq}, we easily get $K_c'(t)=2(2-\beta)f_c(t)f'_c(t)^2$.
Assume now that $f'_c$ vanishes, and let $s_c$ be the first point such that $f'_c(s_c)=0$. Then $f'_c$ and $f_c$ are positive on $[0,s_c)$, and hence $K_c$ is nondecreasing on $[0,s_c]$. Since $f''(s_c)\leq0$, we have $K_c(s_c)=2f_c(t)f''_c(t)-\beta f_c(t)^2<0$. This implies that $K_c(0)<0$. 

Consequently, if $K_c(0)\geq 0$, then $f_c'>0$ on $[0,T_c)$. From Proposition~\ref{paa}, it follows that $T_c=+\infty$ and  $f'_c$ tends to $0$ or $1$ at infinity. But, if $f'_c(t)\to 0$ as $t\to+\infty$, then we obtain a contradiction as above, since $K_c(t)\to -\beta\,l^2$ as $t\to+\infty$, where $l$ is the limit of $f_c$ at infinity (see Lemmas~\ref{ld} and~\ref{le}). The proof is complete, since $K_c(0)=2ac-b^2+(2b-\beta)a^2$.
\ep

\coro{If $\beta\in(1,2]$, $a>0$ and $b>0$, then the problem $({\cal P}_{\beta;a,b,1})$ has infinitely many solutions. }
\pr This follows immediately from Proposition~\ref{pab}.
\ep


\bigskip

Mohammed AIBOUDI 

D\'epartement de Math\'ematiques

Facult\'e des Sciences Exactes et Appliqu\'ees

Universit\'e d'Oran1, Ahmed Benbella

Oran, ALG\'ERIE  

E-mail: {\tt m.aiboudi@yahoo.fr}

\bigskip

Ikram BENSARI n\'ee KHELIL 

D\'epartement de Math\'ematiques

Facult\'e des Sciences Exactes et Appliqu\'ees

Universit\'e d'Oran1, Ahmed Benbella

Oran, ALG\'ERIE 

E-mail: {\tt bensarikhelil@hotmail.fr}

\bigskip

Bernard BRIGHI 

Laboratoire de Math\'ematiques, Informatique et Applications 

Universit\'e de Haute-Alsace

4 rue des fr\`eres Lumi\`ere 

68093 Mulhouse, 
FRANCE 

E-mail:  {\tt bernard.brighi@uha.fr}


\begin{thebibliography}{99}      


\bibitem{a} R.C. Ackerberg, {\em Boundary layer separation at a free stream-line}, J. Fluid 
Mech. \textbf{44} (1970) 211-226.

\bibitem{ab}  {M. A\"iboudi and B. Brighi}, {\em On the solutions of a boundary value problem arising in free convection with prescribed heat flux}, {Arch. Math.} \textbf{93}(2) (2009) 165-174.                                   

\bibitem{aly} {E.H. Aly, L. Elliott and D.B. Ingham}, {\em Mixed convection boundary-layer flows over a vertical surface embedded in a porous medium}, {European J. Mech. B/Fluids} 
\textbf{22} (2003) 529-543.       
                                    
\bibitem {bbt} {Z. Belhachmi, B. Brighi and K. Taous}, {\em On a family of differential equations for boundary layer approximations in porous media}, {European J. Appl. Math.} \textbf{12}(4) (2001) 513-528.

\bibitem{bk} I. Bensari-Khelil, {\em Ph'D Thesis}. In preparation. 

\bibitem{bla} {H. Blasius}, {\em Grenzschichten in Fl\"ussigkeiten mit kleiner Reibung}, {Z. Math. Phys.} \textbf{56} (1908) 1-37.

\bibitem {bb} {B. Brighi}, {\em On a similarity boundary layer equation}, {Z. Anal. Anwendungen} \textbf{21}(4) (2002) 931-948.

\bibitem {asmq} {B. Brighi}, {\em Sur un probl\`eme aux limites associ\'e \`a l'\'equation diff\'erentielle $f'''+ff''+2f'^2=0$}, {Ann. Sci. Math. Qu\'ebec}  \textbf{33}(1) (2009) 23-37.

\bibitem{rm} B. Brighi,  {\em The 
equation $f'''+ff''+{\bf g}(f')=0$ and the associated boundary value problems}, {Results. Math.} \textbf{61}(3-4) (2012) 355-391.                              

\bibitem {bfs} {B. Brighi, A. Fruchard and T. Sari}, {\em On the Blasius problem}, {Adv. Differential Equations} \textbf{13}(5-6) (2008), 509-600. 

\bibitem {mmas} {B. Brighi and J.-D. Hoernel}, {\em On similarity solutions for boundary layer flows with prescribed heat flux}, {Math. Methods Appl. Sci.} \textbf{28}(4) (2005) 479-503.

\bibitem {aml} {B. Brighi and J.-D. Hoernel}, {\em On the concave and convex solutions of mixed convection boundary layer approximation in a porous medium}, {Appl. Math. Lett.} \textbf{19} (2006) 69-74.

\bibitem {amuc} {B. Brighi and J.-D. Hoernel}, {\em On a general similarity boundary layer equation}, {Acta Math. Univ. Comenian.} \textbf{77}(1) (2008) 9-22. 

\bibitem {bs} {B. Brighi and T. Sari}, {\em Blowing-up coordinates for a similarity boundary layer equation},  {Discrete Contin. Dyn. Syst. A} \textbf{12}(5) (2005) 929-948.

\bibitem {bt} {B. Brighi and J.-C. Tsai}, {\em Similarity solutions arising from a model in high frequency excitation of liquid metal with an antisymmetric magnetic field}, IMA J. Appl. Math. \textbf{77} (2012) 157-195.

\bibitem{cheng} {P. Cheng and W. J. Minkowycz}, {\em Free-convection about a vertical flat plate embedded in a porous medium with application to heat transfer from a dike}, {J. Geophys. Res.} \textbf{82}(14) (1977) 2040-2044.

\bibitem{falk} {V. M. Falkner and S. W. Skan}, {\em Solutions of the boundary layer equations}, {Phil. Mag.} \textbf{7}(12) (1931) 865-896.

\bibitem{g} {M. Guedda}, {\em Similarity solutions of differential equations for boundary layer approximations in porous media}, {J. Appl. Math. Phys.  (ZAMP)} \textbf{56} (2005) 749-762.

\bibitem{guedda} {M. Guedda}, {\em Multiple solutions of mixed convection boundary-layer approximations in a porous medium}, {Appl. Math. Lett.} \textbf{19} (2006) 63-68.

\bibitem {hart} {P. Hartman}, {Ordinary Differential Equations}, Wiley, New-York, 1964.
 
\bibitem{kyy} {H.C. Kang, J.C. Yang and G.C. Yang}, {\em Existence and uniqueness of concave and convex solutions of mixed convection equation}, Nonlinear Analysis Forum \textbf{13}(2) (2008) 157-165.

\bibitem{m} {J. B.  McLeod}, {\em The existence and uniqueness of a similarity solution arising from separation at a free stream line}, {Quart. J. Math. Oxford Ser. (2)} \textbf{23} (1972) 63-77.

\bibitem{p} J.E. Paullet, {\em An uncountable number of solutions for a BVP governing Marangoni convection}, Math. Comput. Modelling  \textbf{52} (2010) 1708Ð1715. 

\bibitem{sc} {B.B. Singh and I.M. Chandarki}, {\em On the asymptotic behaviours of solutions of third order non-linear differential equation governing the MHD flow}, {Differ. Equ. Appl.} \textbf{3}(3) (2011) 385-397. 

\bibitem{tw} {J.-C. Tsai and C.-A. Wang}, {\em A note on similarity solutions for boundary layer flows with prescribed heat flux}, {Math. Methods Appl. Sci.} \textbf{30}(12) (2007) 1453-1466.

\bibitem{y} {G.C. Yang}, {\em An extension result of the opposing mixed convection problem arising in boundary layer theory}, Appl. Math. Lett. \textbf{38} (2014) 180-185.

\bibitem{yzd} {G.C. Yang, L. Zhang and L.F. Dang}, {\em Existence and nonexistence of solutions on opposing mixed convection problems in boundary layer theory}, European J. Mech. B/Fluids \textbf{43} (2014) 148-153.

\end{thebibliography}
\end{document}